# Towards a Critical Pragmatic Philosophy of Sustainable Mathematics Education

Dennis Müller

University of Cologne, Germany

**Abstract:** This paper proposes critical pragmatism as a philosophy of sustainable mathematics education to bridge the gap between critical theory and the existing patchwork implementations. Combining existential sustainability as a holistic concept with pragmatic frameworks from the ethics in mathematics education literature creates a foundation enabling critical reflection and pragmatic implementation. We outline how their synthesis naturally leads to a three-stage implementation strategy: cultivating an ethical classroom culture, engaging with ethnomathematics, and tackling complex sustainability problems. Our critical pragmatic approach attempts to build a new philosophical perspective to equip teachers and students with the mathematical competencies, critical perspectives, and ethical grounding necessary to navigate and contribute to a sustainable future and to provide new analytic pathways.

## Introduction: The Critical-Pragmatic Challenge in Sustainable Mathematics Education

We live in the Anthropocene, an era demanding transformative action in response to unforeseen, interconnected socio-ecological crises (Frankopan, 2023), pushing us to see nature beyond an object governed by universal scientific laws (Latour, 2017). Education is not just one of the 17 Sustainable Development Goals (United Nations, 2015); it is also recognized as an indispensable tool for empowering individuals with the knowledge, moral values, and agency needed to achieve a more sustainable world. Education has become both a goal and the means to a sustainable end (Cai & Wolff, 2023). However, even after decades of research into climate change and other sustainability concerns, a persistent societal gap exists between awareness and action (Beckert, 2024), suggesting that education must move beyond mere information transmission to foster civic dialogue and action (Brennan, 2017).

Within mathematics education, the integration of sustainability concerns still faces significant challenges. At present, the implementation of sustainable mathematics education[1] into curricula and classrooms frequently remains superficial and focused on selective skills without deeper critical engagement (Bamber et al., 2016; Lestari et al., 2024; Tesfamicael & Enge, 2024; Wilhelm, 2024). Juxtaposed with these patchwork efforts is the work of critical mathematics scholarship, including specific arguments focused on post-normal science education (Barwell, 2013, 2018; Steffensen, 2017) and eco-justice (Wolfmeyer et al., 2017), challenging the modern enterprise of mathematical science at a fundamental level (Valero & Beccuti, 2025). These perspectives scrutinize institutions, power structures, and unsustainable ways of life more generally while aiming to empower learners for ecological, social, and political transformation. However, evidently, the profound critiques offered by critical mathematics education remain difficult to translate into teaching and teacher training (Su et al., 2023; Vásquez et al., 2023; Vásquez et al., 2020), leading to a gap between critique and action (see also Valero & Beccuti, 2025). To make matters more difficult, sustainable mathematics education itself appears philosophically underdeveloped (Li & Tsai, 2020), often solely focused on critical mathematics education (cf. Makramalla et al., 2025; Skovsmose, 2023, 2024), so that alternative perspectives remain rare.

This dichotomy reveals a core tension resonant with the broader philosophical dialogue between Critical Theory and Pragmatism, here specified to sustainable mathematics education (see also Müller, 2024). How can mathematics education be both critically insightful, addressing systemic issues and power dynamics, and even promoting system change (the "critical" imperative), and pragmatically effective, offering workable pedagogical approaches that address real-world problems and foster agency within existing constraints (the "pragmatic" imperative)?

We argue that critical pragmatism offers a valuable philosophical lens for navigating this tension. Critical pragmatism seeks to synthesize the strengths of both traditions: pragmatism's focus on experience and practically minded action, combined with critical theory's commitment to analyzing power, promoting social justice, and pursuing emancipation. It aims for a grounded philosophy, reflecting both needs and constraints in an ethically responsible, critically aware, and transformation-oriented perspective (Rezac & Husbye, 2024; Ulrich, 2007).

---

[1] In this paper, we will speak of sustainable mathematics education rather than Mathematics Education for Sustainable Development, as MESD reflects a particular kind of sustainability education typically associated with the Sustainable Development Goals and the Agenda 2030. However, we do not want to restrict our approach in this way and consider sustainable mathematics beyond temporarily predetermined views.

We propose a framework integrating holistic perspectives on existential sustainability (Persson, 2024) with pragmatic ideas from the ethical turn in mathematics education (Ernest, 2024; Müller, 2025) that functions as a specific, potent form of critical pragmatism tailored to the challenges of sustainable mathematics education. We will argue that existential sustainability, focusing on the meanings and values of existential concerns within sustainability, provides the critical depth and breadth to cover existing sustainability conceptualizations within mathematics education. Simultaneously, an explicit focus on ethical practice emphasizes individual responsibility, offering the normative grounding and pragmatism required to orient existential sustainability towards a functioning classroom experience. Their synthesis, we contend, acts as a critical pragmatic bridge. It addresses critiques raised by critical scholarship regarding power, mathematics' neutrality, and modernity's unsustainability while providing actionable pathways for implementation.

In the following, we will review existing conceptualizations of sustainable mathematics education before we elaborate on existential sustainability and ethics through a critical pragmatic lens, and finally outline a sequenced implementation strategy designed to enact this synthesis in practice. Beginning with establishing sustainable and ethically-minded classroom cultures (e.g., as conceptualized by Narvaez, 2010; Wilhelm, 2024), ethnomathematical practices are introduced to get students to appreciate diverse mathematical perspectives, cultural contexts and new epistemic foundations (Machaba & Dhlamini, 2021). To do so, ethnomathematics will be understood as "embedded in ethics, focused on the recovery of the cultural dignity of human being" (D'Ambrosio, 2006, p. 1). Complex sustainability problems, such as climate change, can then be holistically explored. Using lessons from existential sustainability and the ethical turn in mathematics education, we will show how these three existing educational methodologies can be pragmatically combined to inform applicable critical teaching of sustainability concerns.

## Existing Conceptualizations of Sustainable Mathematics Education Research

The subject's most prevalent conceptualizations involve the three domains of environmental, social, and economic sustainability (Kabul & Kaleci, 2024), reflecting the interconnected web of the 17 Sustainable Development Goals. In this context, Skovsmose (2023, p. 47) notes that sociopolitical and environmental concerns cannot be separated when one is to educate for sustainable transformation. This tripartite structure is then complemented by Sterling's (2001) taxonomy of accommodation, reformation, and transformation, leading to mathematics education about, for, and as sustainability (Renert, 2011). Because Sterling's and Renert's distinction is connected to first, second, and third-order learning, this approach also appears in other parts of critical mathematics education, such as justice-driven mathematics, where authors also readily speak of mathematics about social justice, for social justice, and as social justice (cf. Stinson et al., 2012; Xenofontos et al., 2021).

Thus, a broad and action-driven perspective on mathematical competencies, for example, as an "insightful readiness to act appropriately in response to a specific sort of mathematical challenge in given situations" (Niss & Højgaard, 2019, p. 14), meets complex critical ideas about sustainability aimed at empowering students with both mathematical competencies and general thinking skills (Oikonomou, 2025; Widiati & Juandi, 2019), including developing critical perspectives of the subject's formatting power in sustainability discourses (Skovsmose, 2023; Steffensen et al., 2021). Building on a long-standing tradition focused on mathematical literacy and responsible citizenship, a key characteristic of critical sustainable mathematics education is the strong emphasis placed on the social dimension, whereby sustainability is generally not viewed as a technocratic management problem solvable by existing modes of modern mathematical reasoning, but as a postmodern educational issue (e.g., Valero & Beccuti, 2025). Such perspectives resist simplistic interpretations and consciously acknowledge that sustainability always requires an integrated approach (Kabul & Kaleci, 2024; Renert, 2011) through an "entangled" perspective (Makramalla et al., 2025).

Existing conceptualizations are then further nuanced by a temporal component. They are generally framed as dynamic minds focused on imagined futures (cf. Beckert, 2016). Aiming at mathematical literacy, the empowerment of students, and societal change, these are then similar to Jeronen's (2013, p. 2371) definition of "sustainability [as] a paradigm for thinking about the future [...] in the pursuit of an improved quality of life" or Pitt's (2009, p. 37) perspective on education for sustainable development as a "frame of mind [...] bringing to the surface underlying values and beliefs through the exploration of contradictions and arguments." General competencies go hand in hand with mathematical competencies to form a dynamic mode of reasoning aimed at future life. Additionally, there is a tendency to think that the problem cannot, or at least will not, be solved using technocratic, top-down approaches and that significant and continuous bottom-up efforts are required:

"One of the lessons of chaos is that creative emergence cannot be controlled top-down. It is a bottom-up project that involves the diverse contributions of many interacting participants. We are these participants - educators, researchers, and students who are passionate about mathematics and the role it can play in the world." (Renert, 2011, p. 25)

All of this means that sustainable mathematics education is as much a vision for the future as it is an internal state full of hope, albeit with the occasional struggles to find it, as it can also happen in other critical scholarship (cf. Ryan & Steffensen, 2021). There is, after all, no action without hope, a concept that can be traced back to the founding figures of critical pedagogy themselves (Gottesman, 2016). Here, sustainable mathematics education does not just share a concern for the internal life with the wider critical mathematics scholarship but also with recent writings on ethics in mathematics (cf. Müller, 2025).

The socio-political perspective (Gutiérrez, 2013) is thus always part of sustainable mathematics education. The ever-present search for the meaning of mathematics education in such political writings (Pais & Valero, 2012) manifests as a "socio-ecological turn" with three strands: (1) interconnected life on earth, (2) the necessity to reevaluate human and non-human relationships on this planet, and (3) the need to redefine the aims of mathematics education through the lens of the first two (Coles, 2022) using (typically) constructivist pedagogical theories focused on real-world problems situated in the students' socio-cultural context.

Globally, but particularly within the Global South and where Indigenous communities are present, we are also seeing an increasing focus on decolonization, Indigenous knowledge systems, localized cultural relevance, and diversified epistemological approaches in both empirical and theoretical scholarship (e.g., Garcia-Olp et al., 2022; Kulago et al., 2021; Nicol et al., 2020; Okeke & Salami, 2017; Stavrou & Miller, 2017). Such approaches may explore the "dualisms [that] separate mind from body, body from nature, and spirit from matter" (Kulago et al., 2021, p. 345) and thus can provide radically different conceptions from traditional mathematics education. These approaches can be difficult to incorporate into classical classrooms, as they require "productive struggles that [require] enhanced higher-order levels of thinking [, learning, and teaching]" (Gula & Jojo, 2024, p. 404). Decolonial theory and teaching practices are known to challenge existing notions of social justice (Fúnez-Flores et al., 2024), and they accomplish a similar task in the mathematical discourse on sustainability.

That these profound critical theories exist alongside a patchwork implementation shows that we are only beginning to "rethink [and] re-envision mathematics teaching and learning for 21st-century learning priorities" (Li & Tsai, 2020, p. 139) and that mathematics education may not yet have found the right conceptual frameworks that simultaneously permit critical analysis, pragmatic teacher training, and implementation strategies within socio-economic and political constraints. Next, we explore how "existential sustainability" might be the first part of a solution.

## Foundation 1: Existential Sustainability

The diversity of scholarship on sustainable mathematics education is connected to varying understandings surrounding existential risks, particularly those concerning climate change, which prompt notably different interpretations across different actors in science, politics, and society. Huggel et al. (2022) demonstrate that the existential risk studies literature, the Intergovernmental Panel on Climate Change (IPCC), and societal climate movements understand such high-end dangers differently. What are these differences? The existential risk literature focuses on the extinction of humanity and intelligent life from natural and anthropogenic causes (Bostrom, 2002). Climate change movements readily engage with such narratives but also come with fundamental justice concerns, capitalist critiques, and calls for system change (Han & Ahn, 2020). The IPCC acts more cautiously in its reports, focusing more on "Reasons for Concern" and limiting existential risks to smaller scales (e.g., Small Island States) (Huggel et al., 2022, pp. 8–9). In short, existential risks are portrayed using various narratives, even when there is a scientific consensus on the existence of climate change itself.

Emerging from sustainability studies, the concept of "existential sustainability" seeks to capture the complexity and heterogeneity within current sustainability discourses (Persson, 2024), doing justice to the multifaceted dimensions of "existential" and "sustainable." Persson (2024, 10) notes its function in "linking sustainability thinking with overlapping discourses in risk research, threat theories, and loss and damage" studies. This idea encompasses varying scopes (individual to planetary), qualities (well-being to biological death), species (human and non-human), and forms of materiality (including life, meaning, values, and cultures). This connection manifests across different levels: psychologically, it addresses existential aspects like freedom, belonging, and anxiety alongside sustainable considerations of long-term personal growth. At the community and cultural level, concerns for values, identity, and connections are balanced against the sustainable goals of intergenerational continuity, resilience, and social cohesion, including a critical stance against homogenization. For the climate and biosphere, the concept considers the interconnectedness of life through an evolutionary lens while promoting the maintenance of populations, ecological resilience, and genetic diversity. Existential sustainability posits an intrinsic right to existence for all forms of life, being, and materiality. It thereby offers an affirmative, multi-perspectival outlook on existence, acknowledging complexities not fully captured by purely Western rational arguments.

Hence, existential sustainability is critical at its very core. It directly challenges the dominance of purposive-rational action in political discourses that attempt to over-rationalize sustainability. Providing material and immaterial beings with independent rights to existence, focused on dignity, values, and meaning, reasserts Habermas's focus on the "interaction" sphere, countering technocratic images of sustainable mathematics education that obscure deeper ethical questions. For example, the survival of different cultures readily opens up existing critiques of power and justice, enabling different epistemic outlooks on life. In doing so, existential sustainability undermines the ideological pretenses of science within capitalist societies (cf. Habermas, 1968), forcing us to reconsider how much mathematics can capture and where its formatting power ends.

However, existential sustainability is also pragmatic. Grounded in localizable knowledge and experience, it takes seriously the fallibility of knowledge systems to reject uncontextualized sustainability claims. Its focus on multidimensional real-world consequences (e.g., the loss of meaning, eco-anxiety, habitat loss, etc.) answers Dewey's (1997, p. 17) warning of extremist "Either-Ors" by providing a breadth that enables students and teachers to pragmatically discover their existential challenges and outlook on sustainability inside a common terminology. It helps to overcome challenges related to classical sustainability concerns that can make it difficult for people to feel connected when their own experiences are not necessarily reflected in these.

However, we cannot stop here just because existential sustainability seems to capture the breadth and depth of our existing discourse. Existential sustainability is, in a very real sense, underdetermined, allowing us to fill it with different interpretations of the "existential" and "sustainable." This flexibility, while valuable, might leave educators - particularly those accustomed to more defined mathematical frameworks - seeking additional normative guidance. To achieve what Makramalla et al. (2025) call a "dialectic between the individual and multiple levels of the social and ecological," we arguably still need something more, and this is where ethics in mathematics education comes into play.

## Foundation 2: Ethics in Mathematics Education

The critical pragmatic paradigm requires value-laden classroom inquiry. While this is already partially encoded in existential sustainability, much of it only emerges in specific student-teacher interactions. Our aim here is not to develop a new normative perspective but to edge closer to a foundation that helps students and educators pragmatically explore and develop their critical worldview. As the emerging ethical turn extends previous sociopolitical writings on mathematics education by sharpening the focus on individual responsibility (Ernest, 2024; Müller, 2025), potentially even positioning ethics as a "first philosophy" for mathematics education (Ernest, 2012), its pragmatic scholarship may offer us a pathway to operationalize existential sustainability.

The scholarship on ethics in mathematics rejects persistent myths of mathematical neutrality, purity, and universality (NPU) (Ernest, 2020b; Müller, 2025), acknowledging instead that all mathematical artifacts have politics (cf. Winner, 1980), requiring the ethical awareness of its practitioners, students, and teachers. At its core, this awareness is not restricted to critical perspectives but can readily include other versions of philosophical ethics. Crucially, it begins with the understanding that mathematics teachers can act at different levels (Rycroft-Smith et al., 2024, p. 361):

- Level -1: Actively obstructing efforts to address ethics in mathematics
- Level 0: Believing there is no ethics in mathematics
- Level 1: Realising there are ethical issues inherent in mathematics
- Level 2: Doing something: speaking out to other mathematics [teachers]
- Level 3: Taking a seat at the tables of power
- Level 4: Calling out the bad mathematics of others

By showing teachers that incremental steps are possible, this level-based understanding is structured around reflecting one's ethical stance and gaining the confidence to discuss ethical and sustainable issues. The latter, in particular, can still be lacking among teachers, as Coles (2023) discusses. It provides a pragmatic understanding of Monroe et al.'s (2019) observation that sustainable education needs deliberate discussion to turn awareness into action. Higher levels of ethical awareness require educators to confront the idea that their teaching, and mathematics education more generally, can harm both students and wider society (Abtahi, 2022; Ernest, 2016, 2018, 2020a, 2020b; Mohanty & Swain, 2018; Müller, 2025; Rycroft-Smith et al., 2024; Schoenfeld, 1988).

In this context, Chiodo and Bursill-Hall (2019) argue that teaching ethics in mathematics in ways that students appreciate, leads to agency, and normalizes ethical considerations, requires separate ethics courses, the interweaving of ethical content into mathematics sessions, and the support of other faculty members. While one cannot introduce a separate course at the school level, all three ideas can still be transferred. The separate course can become a project week, the interweaving can involve regularly selecting mathematical problems with a sustainability component, and the support of the faculty turns into the support of other teachers. This approach towards embedded forms of ethical teaching is pragmatic compared to other attempts, and the interweaving of ethical content mirrors a theme that is also present in sustainable mathematics education, where a "one-off" feel can develop if topics are not sufficiently integrated (Coles, 2023). While the levels of ethical awareness are centered on the teacher, Chiodo and Bursill-Hall (2019) provide a wider structure on where and when the engagement happens. What is left is finding a pragmatic perspective for how this engagement can happen. The literature has proposed different normative mathematical lifecycles (e.g., Pohlkamp & Heitzer, 2021). However, as we are aiming to connect this foundation to existential sustainability later, we focus on a practically informed "pillar-based approach" that was also developed in the context of existential risks (Chiodo & Müller, 2023):

1. Deciding whether to begin
2. Diversity and perspectives
3. Handling data and information
4. Data manipulation and inference
5. The mathematisation of the problem
6. Communicating and documenting your work
7. Falsifiability and feedback loops
8. Explainable and safe mathematics
9. Mathematical artefacts have politics
10. Emergency response strategies

What these pillars show us is threefold: (1) ethical questions are indeed part of the entire mathematical workflow, (2) the early and late stages of a mathematical work can be less mathematical and more political, and (3) moving away from the mathematics (measured by distance to pillar five) potentially requires both higher levels of ethical awareness and better institutional support.

These ethics frameworks provide a pragmatic perspective regarding teachers, institutional settings, and mathematical workflows. They can help educators understand their ethical obligations as humans, citizens, and mathematics teachers (Ernest, 2019a, 2019b), including exploring the duality of "doing good" and "preventing harm" (Müller, 2025). Additionally, the pragmatic literature shows that ethical engagement can essentially happen with three interconnected foci: mathematics, its community of practitioners (e.g., students), and wider societal and planetary concerns (Chiodo & Müller, 2024). Such a range of concerns also aligns with the multidimensionality of existential sustainability. In the following, these works will guide us to find a suitable place between what Gutiérrez (2009) described as "playing the game" and "changing the game" of mathematics.

## Implementing Sustainable Mathematics Education

Having established the foundations for a critical pragmatic philosophy of sustainable mathematics education, we now transition from theoretical exploration to practical implementation concerns, outlining a sequenced pedagogical strategy. We attempt to offer a pragmatic framework for addressing the complex, heterogeneous, and potentially overwhelming concerns encompassed by existential sustainability. To maintain a critical dimension, we use sustainable classrooms, focused on students' well-being and gentle teaching to represent community concerns, ethnomathematics as an example concerning the ethical development of mathematics, and mathematical climate change models for wider planetary and societal concerns. Each of these represents different aspects of existential sustainability. The sustainable classroom is immediately connected to a wholesome picture of students, ethnomathematics relates to the sustainable development of mathematics and the preservation of different cultures, and climate change education is about risks in a more classical sense. We present these in a sequenced approach - (1) classroom culture, (2) ethnomathematics, (3) planetary/societal sustainability concerns - and argue that their progression can be understood as a pathway intended to move students and educators beyond lower levels of ethical awareness.

In practice, it is likely best to consider interleaving the cultural and epistemic insights from ethnomathematics (stage 2) with the critical modeling of complex sustainability issues (stage 3) rather than treating them as separate steps. However, we will consider them as properly sequenced for the rest of this article for simplicity and to provide a cleaner argument. This deliberate sequencing of stages reflects a core pedagogical rationale whereby robust affective grounding within the immediate community of learners (stage 1) is indispensable for navigating sensitive cultural concepts (stage 2) and tackling the ethically charged nature of complex goals such as the mitigation of climate change (stage 3). Climate change, in particular, is an issue that includes physical (e.g., harm from high temperatures), emotional (e.g., eco-anxiety), and immaterial (e.g., the loss of one's culture) existential sustainability risks for individual students. Even if we only consider eco-anxiety and feelings of helplessness, we can see why having a deeper grounding in communal and cultural ethics may be necessary first. Not pathologizing eco-anxiety and building emotional resilience is essential for young people to foster constructive engagement (Brophy et al., 2023). However, this can only happen inside the classroom if its culture is already sustainable and provides the space for such emotions, dictating that the sustainable classroom comes before attempts to teach complex sustainability issues.

By fostering their well-being and sense of belonging, this first stage enables a search for meaning by directly engaging with existential sustainability's psychological and interpersonal dimensions. Additionally, this is likely the first contact for students (and maybe even teachers) that actively challenges the level 0 perspective that mathematics and its teaching is a value-free, ethically pure endeavor by making it an explicit topic inside the classroom. While not directly focused on the mathematical workflow's technical (middle) pillars, this stage focuses on social and socio-mathematical norms to build interpersonal trust and ethical sensitivity necessary for explorations beyond the classroom itself. Reducing anxiety and increasing empathy with an openly communicated commitment to psychological safety is essential for students to develop the competencies necessary for more complex explorations and agency. In this context, Wiemers et al. (2012) stress the benefit of collegiate dialogue (similar to Rycroft-Smith et al.'s (2024) level 3), opening the perspective for (critical) theories of communicative action and discourse ethics.

Research on the social and socio-mathematical norms in mathematics education shows that students and teachers come together in complex interactions to find these norms through communication and action (Meyer & Schwarzkopf, 2025). Pragmatically, this involves establishing classroom conditions and communicative practices that approximate an ideal inclusive speech situation free from coercion, where learners and teachers have the opportunities to express themselves freely. Explicitly embedding the creation and justification of norms within the classroom practice enables students and teachers to move from level 0 or level 1 to higher levels.

In short, stage 1 begins to operationalize "mathematics education as sustainability." Pragmatic mathematical activities focused on the student's needs and desires are then easier to explore, making space for real-world mathematics exercises that the literature so often calls for. One concrete vision for this can be found in Wilhelm (2024).

Regarding stage 2, Borba (1990) identifies the "efficiency of ethnomathematics" in its cultural roots while warning that students may not necessarily be interested in learning the ethnomathematics of other cultures (ibid., p. 41), and more recently, Wulandari et al. (2024) surveyed ethnomathematical research, finding sufficient evidence that mathematics rooted in the students' culture is more appreciated than contextless mathematics. However, this is a problem for Western classrooms as their ethnomathematics is the "Western way of mathematics," which is typically understood as neutral, pure, and universal (Ernest, 2020b; Müller, 2025). Nonetheless, the "socio-critical modeling" (Kaiser & Sriraman, 2006, p. 306) that is necessary for sustainable mathematics education traces back to ethnomathematical perspectives of mathematics (Coles, 2016, p. 3). Thus, the question is no longer whether to introduce ethnomathematics as the next step but how to do it.

A critical pragmatic implementation avoids seeing ethnomathematics merely as a means for a higher end and always as an end in itself. Thus, following Chiodo and Bursill-Halls's (2019) tripartite approach to ethical education (projects, regular exercises, and teacher support), we suggest an independent project on ethnomathematics, the regular selection of ethnomathematical exercises, and for teachers to engage at levels 2 and 3 of ethical awareness, i.e., to speak about different mathematical cultures within, and outside, the classroom. Leveraging ethnomathematics this way is aimed at cultivating critical thinking skills beyond those of one's own culture and potentially paves the path for the epistemic humility needed for more systematic societal change.

"Ethnomathematics makes it clear that mathematics and mathematical reasoning are cultural constructions. This raises the challenge to embrace the global variety of cultures of mathematical activity and to confront the politics that would be unleashed by such attention in a typical North American school. That is, ethnomathematics demands most clearly that critical thinking in a mathematics classroom is a seriously political act" (Appelbaum, 2004; as cited in François & Stathopoulou, 2012, p. 237).

While ethnomathematics centers on pillars 2 and 9, this political act is one prerequisite for a deeper engagement with socio-critical modeling (as represented by pillars 3-8). Socio-critical modeling represents a critically informed engagement with the third core area of ethical concern - mathematics's role in broader planetary and societal issues - while simultaneously confronting existential sustainability concerns at their widest scope (e.g., climate change, biodiversity loss). Having developed emotional resilience and a sustainable classroom (stage 1) and significantly broadened their cultural and epistemic perspectives (stage 2), students are now positioned to deploy their mathematical competencies critically and meaningfully to real-world issues. Addressing the interconnected sustainability problems that the world is facing necessitates holistic scientific competencies within authentic, challenging contexts (Reiss, 2023). As for ethnomathematics, we suggest that project-based learning of sustainability must be complemented by interwoven sustainability exercises outside of projects and within normal classroom situations, as only this can normalize ethical and, thus, sustainable thought in ways that foster student agency beyond the classroom.

Drawing upon socio-critical modeling perspectives, students learn to "unbox mathematics" and develop a deeper awareness of the formatting power of the subject (Gibbs, 2019, p. 16). This critical stance requires direct engagement with multiple ethical pillars: careful consideration of data handling (Pillar 3) and inference (Pillar 4), deep reflection on the act of mathematization itself (Pillar 5), ensuring clear communication (Pillar 6), maintaining openness to scrutiny (Pillar 7), promoting explainability (Pillar 8), and explicitly recognizing that models possess political dimensions (Pillar 9). The ability to interrogate a model's underlying assumptions, data, values, biases, and societal impact relies on the epistemic humility and appreciation for diverse perspectives that were cultivated in earlier stages, as well as the classroom culture of stage 1, without which students would likely not fully and honestly take on this challenge. Successfully navigating this stage provides the context and analytical tools necessary for higher ethical awareness and action (e.g., level 4) (cf. Rycroft-Smith et al., 2024).

In summary, we propose to jointly explore sustainable classrooms, ethnomathematics, and classical sustainability through pragmatic means:

1. Stage: Cultivating sustainable and ethical classroom cultures.
    a. Primary ethical locus: Community of mathematical practitioners (students, teachers).
    b. Pragmatic goal: Addressing ethical responsibilities within the immediate community and attending to their existential sustainability concerns.
    c. Level of awareness: Challenging Level 0 and establishing Level 1 within the community to build the foundations for more collaborative ethical practice.
2. Stage: engaging with ethnomathematics for cultural sustainability.
    a. Primary ethical locus: The discipline of mathematics itself.
    b. Pragmatic goal: Addressing ethical questions concerning the nature and development of mathematics and promoting the existential sustainability of diverse knowledge systems.
    c. Levels of ethical awareness: A deepened awareness of level 1 fosters the critique necessary for higher levels; it directly engages with Pillar 2: Diversity and perspectives and explores Pillar 9: Artefacts have politics.
3. Stage: Tackling classical sustainability concerns.
    a. Primary ethical locus: Planetary and societal sustainability concerns.
    b. Pragmatic goal: Addressing humanity's impact on wider planetary and societal issues, exploring the formatting power of mathematics and its social justice implications, and confronting large-scale risks.

c. Levels of ethical awareness: Providing context and tools for a deeper engagement with levels 2, 3, and 4; deeper engagement with pillars 3-9, particularly critical mathematization, modeling, and evaluation.

## Closing Remarks

This paper argued for a critical pragmatic philosophy to navigate the complexities of sustainable mathematics education in the era of the Anthropocene. Recognizing the gap between awareness and action in addressing socio-ecological crises and the limitations of current approaches in mathematics education, we proposed a synthesis of two core foundations: existential sustainability and the ethical turn in mathematics education.

Existential sustainability offers a multidimensional lens, encompassing diverse concerns from individual well-being and cultural dignity to planetary boundaries. Stories are key in sustainability education (Amico et al., 2023), and as we have argued, existential sustainability allows us to build stories that create an intimate relationship between students' existential questions and larger sustainability discourses. Existential sustainability acknowledges the varying scopes, qualities, and materialities of existence threatened by unsustainability, and it challenges rationalist and purely instrumental views by valuing intrinsic worth, including those of the students. Complementing this, pragmatic scholarship from the ethical turn in mathematics education provides the necessary practical orientation. To create this orientation, we selected scholarship that attempts to be (somewhat) politically neutral. We did so, hoping that this pragmatic perspective enables teachers and their students to find their way. We wanted to provide structure and orientation for asking questions and finding political and ethical answers rather than providing these through a more normative stance.

Our proposed three-stage implementation strategy – (1) cultivating an ethical classroom culture, (2) engaging with ethnomathematics for cultural sustainability, and (3) tackling complex problems like climate change through critical modeling - provides our proposed roadmap. We argued that this sequenced approach fosters affection, empathy, epistemic humility, and critical socio-mathematical competencies in teachers and learners. We contend that integrating these foundations creates a robust, critically informed, pragmatic perspective for sustainable mathematics education by having implicitly argued that together, they are easier to implement than when considered independently.

Critical theory has had a long trajectory from its Marxist origins, incorporating a paradigm shift focused on critiquing (capitalist) civilization under the likes of Adorno and Marcuse and its consequent close connection to protest movements (Schwandt, 2010). Critical pedagogy has a similarly varied tradition (Gottesman, 2016). The critical pragmatic approach taken in this paper was, in some parts, inspired by a personal difficult trajectory of learning about critical theory and critical mathematics education while being confronted with having to teach about sustainability and ethics after having learned mathematics through a traditional curriculum. It was done from the perspective that the critical scholarship questions the pragmatist, and the pragmatic scholarship questions the critique, (both) hoping to find a more sustainable mathematics.